\documentclass[a4paper,12pt]{article}
\usepackage[english]{babel}
\usepackage[T1]{fontenc}
\usepackage{amsthm,enumerate,amsmath,amscd,amssymb}
\usepackage{graphicx}

\newcommand{\R}{\mathbb{R}}
\newcommand{\Rn}{{\mathbb{R}^n}}
\newcommand{\Hn}{{\mathbb{H}^n}}

\newcommand{\arcosh}{\textnormal{arcosh} \,}

\newcommand{\PS}{ {\Rn \setminus \{ 0 \}} }

\newcommand{\comment}[1]{}

\newcounter{minutes}
\divide\time by 60
\newcounter{hours}
\multiply\time by 60
\addtocounter{minutes}{-\time}

\swapnumbers
\numberwithin{equation}{section}
\theoremstyle{plain}

\newtheorem{theorem}[equation]{Theorem}
\newtheorem{corollary}[equation]{Corollary}
\newtheorem{proposition}[equation]{Proposition}
\newtheorem{lemma}[equation]{Lemma}

\newtheorem{remark}[equation]{Remark}

\begin{document}

\begin{center}
{\bf \large Inclusion relations of hyperbolic type metric balls}
\end{center}

\vspace{1mm}

\begin{center}
{\large Riku Klén and Matti Vuorinen}
\end{center}

\vspace{1mm}

\begin{abstract}
  We will consider inclusion of metric balls defined by the quasihyperbolic, the $j$-metric and the chordal metric. The quasihyperbolic metric and the $j$-metric are considered in general subdomains of $\mathbb{R}^n$ and in some particular domains like $\mathbb{R}^n \setminus \{ 0 \}$ and the upper half-space.

  \hspace{5mm}


\end{abstract}

\begin{center}
\texttt{File:~\jobname .tex,
          printed: \number\year-\number\month-\number\day,
          \thehours.\ifnum\theminutes<10{0}\fi\theminutes}
\end{center}

\section{Introduction}

In classical complex analysis several metrics are used in addition to the Euclidean metric. Perhaps the most important are the hyperbolic and the chordal metrics \cite{hmm,kl}. While studying quasiconformal mappings in $\mathbb{R}^n, n \ge 2,$ F.W. Gehring and B.P. Palka \cite{gp} introduced the quasihyperbolic metric which has become one of the standard metrics in higher dimensional geometric function theory \cite{h,vq,vu3}. Due to the absense of a useful counterpart of the Riemann mapping theorem in the higher dimensional context, the hyperbolic metric is not sufficient for the needs of mapping theory, and consequently the quasihyperbolic and related ''hyperbolic type'' metrics (such as the distance ratio and the Apollonian metrics) have been studied by many authors \cite{himps}. General properties of metric structures have been studied recently by A. Papadopoulos and M. Troyanov in several papers, see e.g. \cite{pt}.

For the formulation of our basic problem we consider metric spaces $(X, d_j)$, $j=1,2,$ and suppose that the metrics are \emph{topologically
compatible}. This means that for a fixed $x\in X $ and $r>0$ there exist $m , M >0$ such that
\[
B_{d_1}(x,m)  \subset B_{d_2}(x,r)  \subset B_{d_1}(x,M)
\]
where $B_{d_1}(x,m)$ stands for the ball $\{ z \in X \colon  d_1(x,z)<m  \}$. Note that usually $m$ and $M$ depend on both $x$ and $r$. The purpose of this paper is to study this type of set theoretic inclusion properties for metric balls for the case when $X$ is a subdomain of $\mathbb{R}^n$ and the metrics are of ''hyperbolic type''.

We say that a curve $\gamma \colon [0,1] \to X$ is a \emph{geodesic segment} of the metric space $(X,d)$, if for all $t \in (0,1)$ we have
\[
  d(\gamma(0),\gamma(t)) + d(\gamma(t),\gamma(1)) = d(\gamma(0),\gamma(1)).
\]
A metric space is \emph{geodesic}, if for every pair of points there exists a geodesic segment joining them.

We now proceed to define several of these metrics which enables us to formulate the main theorem.

For a domain $G \subsetneq \Rn$, $n \ge 2$ and a continuous function $w \colon G \to (0,\infty)$ we define the \emph{$w$-length} of a rectifiable arc $\gamma \subset G$ by
\[
  \ell_w(\gamma) = \int_{\gamma}w(z)|dz|,
\]
and the \emph{$w$-metric} by
\begin{equation}\label{wmetric}
  m_w(x,y) = \inf_\gamma \ell_k(\gamma),
\end{equation}
where the infimum is taken over all rectifiable curves in $G$ joining $x$ and $y$.

The \emph{quasihyperbolic metric} is obtained from \eqref{wmetric} with $w(x) = 1/d(x)$, where $d(x)$ is the Euclidean distance between $x$ and $\partial G$, and we denote this metric by $k_G$. By \cite{go} the metric space $(G,k_G)$ is always \emph{geodesic}. Note that in a half-space the quasihyperbolic metric coincides with the hyperbolic metric and therefore the formula of the hyperbolic metric \eqref{hypmetric} apply to the quasihyperbolic metric as well.

The \emph{distance ratio metric} or \emph{$j$-metric} in a proper subdomain $G$ of the Euclidean space $\Rn$, $n \ge 2$, is defined by
$$
  j_G(x,y) = \log \left( 1+\frac{|x-y|}{\min \{ d(x),d(y) \}}\right),
$$
where $d(x)$ is the Euclidean distance between $x$ and $\partial G$.

If the domain $G$ is understood from the context we use the notation $j$ instead of $j_G$ and $k$ instead of $k_G$. The distance ratio metric was first introduced by F.W. Gehring and B.G. Osgood \cite{go} and in the above form by M. Vuorinen \cite{vu2}. The metric space $(G,j_G)$ is not geodesic for any domain $G$ \cite[Theorem 2.10]{k1}.

The \emph{chordal metric} in $\overline{\Rn} = \Rn \cup \{ \infty \}$ is defined by
\[
  q(x,y) = \left\{ \begin{array}{ll}
    \displaystyle \frac{|x-y|}{\sqrt{1+|x|^2}\sqrt{1+|y|^2}}, & x \neq \infty \neq y,\\
    \displaystyle \frac{1}{\sqrt{1+|x|^2}}, & y=\infty.
  \end{array} \right.
\]

The following theorem is our main result and it is proved at the end of Section \ref{halfspace}.
\begin{theorem}\label{mainthm}
  Let $G$ be a proper subdomain of $\Rn$, $x \in G$ and denote $r_q(x) = \min \{ 1/\sqrt{1+|x|^2},|x|/\sqrt{1+|x|^2} \}$.
  \begin{itemize}
  \item[1)] Let $G = \Rn \setminus \{ 0 \}$ and $r \in (0,\pi/2)$. Then
  \[
    B_j(x,m_1) \subset B_k(x,r),
  \]
  where
  \[
    m_1 = \log \left( 1+2 \sin \frac{r}{2} \right).
  \]

  \item[2)] For $G = \Rn \setminus \{ 0 \}$ and $r \in (0,r_q(x))$ we have
  \[
    B_j(x,m_2) \subset B_q(x,r) \quad \textrm{and} \quad B_k(x,m_2) \subset B_q(x,r),
  \]
  where
  \[
    m_2 = \log \left( 1+\frac{2r^2}{\sqrt{1-r^2}} \right).
  \]

  \item[3)] Let $G = \Hn$ and $r > 0$. Then
  \[
    B_j(x,m_3) \subset B_k(x,r),
  \]
  where
  \[
    m_3 = \log \left( 1+\sqrt{2} \sqrt{\cosh r-1} \right).
  \]

  \item[4)] For $G = \Hn$, $r \in (0,r_q(x))$ and $x$ with $x_1 = x_2 = \cdots = x_{n-1} = 0$ we have
  \[
    B_j(x,m_4) \subset B_q(x,r) \quad \textrm{and} \quad B_k(x,m_5) \subset B_q(x,r),
  \]
  where
  \[
    m_4 = \log \left( 1+ \frac{2r}{1-r^2} \right), \quad m_5 = \log \left( 1+\frac{2r^2}{\sqrt{1-r^2}} \right).
  \]
  \end{itemize}
\end{theorem}

It is obvious that this type of theorem is expected to hold not only for the metrics considered above but for numerous other metrics as well. One step in this direction is \cite{m}. We hope to return to this topic in later papers.

\section{General domain}

Throughout the paper $B^n(x,r)$ denotes the Euclidean open ball in $\Rn$.

In this section we consider inclusions of the quasihyperbolic and the $j$-metric balls in general subdomains of $\Rn$. F.W. Gehring and B.P. Palka have showed \cite[Lemma 2.1]{gp} that for any domain $G \subsetneq \Rn$
\[
  j(x,y) \le k(x,y)
\]
for all $x,y \in G$, which is equivalent to
\begin{equation}\label{gehringpalka}
  B_k(x,r) \subset B_j(x,r)
\end{equation}
for all $x \in G$ and $r>0$. The following proposition is an easy consequence of \cite[Theorem 3.8]{s} and \cite[(3.9)]{vu3}.

\begin{proposition}\label{jkjkgen}
  Let $G \subsetneq \Rn$ be a domain and $r \in (0,\log 2)$. Then
  \[
    B_j(x,m) \subset B_k(x,r) \subset B_j(x,r) \subset B_k(x,M)
  \]
  where
  \[
    m = \log (2-e^{-r})
  \]
  and
  \[
    M = \log \frac{1}{2-e^r}.
  \]
  Moreover, the second inclusion is sharp and $M/m \to 1$ as $r \to 0$.
\end{proposition}
\begin{proof}
  Clearly $m > 0$ and since $r \in (0,\log 2)$ we also have $M > 0$.
  By \cite[Theorem 3.8]{s}
\begin{equation}\label{gen2}
B_j(x,r) \subset B^n(x,(e^r-1)d(x))
\end{equation}
  and by \cite[(3.9)]{vu3}
\begin{equation}\label{gen3}
    B^n(x,(1-e^{-r})d(x)) \subset B_k(x,r) \subset B^n(x,(e^r-1)d(x)).
\end{equation}
  By \eqref{gehringpalka}, \eqref{gen2} and \eqref{gen3} we have
  \begin{eqnarray*}
    B_j(x,m) & \subset & B^n(x,(e^m-1)d(x)) = B^n(x,(1-e^{-r})d(x)) \subset B_k(x,r)\\
    & \subset & B_j(x,r) \subset B^n(x,(e^r-1)d(x)) = B^n(x,(1-e^{-M})d(x))\\
    & \subset & B_k(x,M).
  \end{eqnarray*}

  Fix $x \in G$ and let $z \in \partial G$ be any point with $d(x) = |x-z|$. Choose $y \in \partial B_k(x,r) \cap [x,z]$. Now $d(y) = d(x)- |x-y|$ and $j(x,y) = k(x,y)$ implying $y \in \partial B_j(x,r)$. Thus, the inclusion $B_k(x,r) \subset B_j(x,r)$ is sharp.

  By l'H\^opital's rule, the assertion follows.
\end{proof}

Proposition \ref{jkjkgen} is true in all subdomains of $\Rn$. In the following sections we consider similar inclusions and find better radii of the metric balls in some specific domains.

\begin{remark}
  Let $m$ and $M$ be as in Proposition \ref{jkjkgen}. We indicate that $m$ and $M$ can be estimated by a linear function. It is possible to show that for $r \in (0,\log 2)$
  \[
    1+\frac{r}{2} \le \frac{r}{m} \le 1+2r
  \]
  and for $r \in (0,1/2)$
  \[
    1+r \le \frac{M}{r} \le 1+3r, \quad \textrm{and} \quad 1+2r \le \frac{M}{m} \le 1+5r.
  \]
\end{remark}

\section{Punctured space}

In this section we will consider inclusion of metric balls defined by the metrics $q$, $j_G$ and $k_G$ for $G = \PS$. Before considering inclusion of metric balls we introduce a lemma for the chordal metric.

\begin{lemma}\cite[Lemma 7.16]{avv}\label{qball}
  For $x \in \Rn$ and $r \in (0,1/\sqrt{1+|x|^2})$ we have
  \[
    B_q(x,r) = B^n(y,s)
  \]
  and for $x \in \Rn$ and $r \in (1/\sqrt{1+|x|^2},1)$
  \[
    B_q(x,r) = \Rn \setminus \overline{ B^n(y,s) },
  \]
  where
  \[
    y = \frac{x}{1-r^2(1+|x|^2)} \quad \textrm{and} \quad s = \frac{r(1+|x|^2)\sqrt{1-r^2}}{1-r^2(1+|x|^2)}.
  \]
\end{lemma}

The chordal ball $B_q(x,r)$ is either a half-space or the complement of a closed Euclidean ball whenever $r \ge 1/\sqrt{1+|x|^2}$. On the other hand, the chordal ball $B_q(x,r)$ contains 0 whenever $r > |x|/\sqrt{1+|x|^2}$. Since we do not want either of the cases to occur it is natural to assume that $r < r_q(x)$ for $r_q(x) = \min \{ 1/\sqrt{1+|x|^2},|x|/\sqrt{1+|x|^2} \}$. Note that the selection of $r$ implies that $r < 1/\sqrt{2}$ and
\begin{equation}\label{boundsforx}
  \frac{r}{\sqrt{1-r^2}} < |x| < \frac{\sqrt{1-r^2}}{r}.
\end{equation}

We will now consider inclusion of the $j$-metric and the quasihyperbolic metric balls.

\begin{theorem}\label{jkj}
  Let $G = \Rn \setminus \{ 0 \}$, $x \in G$ and $r \in (0,\pi/2)$. Then
  \[
    B_j(x,m) \subset B_k(x,r),
  \]
  where
  \[
    m = \log \left( 1+2 \sin \frac{r}{2} \right).
  \]
  Moreover, the inclusion is sharp and $r/m \to 1$ as $r \to 0$.
\end{theorem}
\begin{proof}
  Let $y \in B_j(x,m)$ and assume first that $|x| \le |y|$. Now $|x-y| < |x|(e^m-1) = 2|x| \sin (r/2) = k(x,z)$ for $z \in \partial B_k(x,r)$ and $|z|=|x|$. For all $u \in \partial B_k(x,r)$ with $|u|>|x|$ we have $|x-u|>|x-z|$ by \cite[Lemma 4.9]{k2} and therefore $y \in B_k(x,r)$. If $|x| > |y|$ we have $y \in B_k(x,r)$ by the previous case and the fact that the quasihyperbolic and $j$-metric balls are invariant under inversion about origin.

  We show that $m$ is sharp. We choose $y \in \partial B_k(x,r)$ with $|y| = |x|$. Now
  \[
    j(x,y) = \log \left( 1+ \frac{|x-y|}{|x|} \right) = \log \left( 1+ \frac{2|x|\sin(r/2)}{|x|} \right) = m
  \]
  and $m$ is sharp.

  Finally, we show that $r/m \to 1$ as $r \to 0$. By l'H\^opital's rule we have
  \[
    \lim_{r \to 0} \frac{r}{m} = \lim_{r \to 0} \frac{1+2\sin(r/2)}{\cos (r/2)} = 1
  \]
  and the assertion follows.
\end{proof}

\begin{corollary}\label{kjk}
  Let $G = \Rn \setminus \{ 0 \}$, $x \in G$ and $r \in (0,\log 3)$. Then
  \[
    B_j(x,r) \subset B_k(x,M),
  \]
  where
  \[
    M = 2 \arcsin \frac{e^r-1}{2}.
  \]
  Moreover, the inclusion is sharp and $M/r \to 1$ as $r \to 0$.
\end{corollary}

\begin{remark}
  Let $m$ be as in Theorem \ref{jkj} and $M$ be as in Corollary \ref{kjk}. It is possible to show that for $r \in (0,\pi/2)$
  \[
    1+\frac{r}{3} \le \frac{r}{m} \le 1+\frac{r}{2}
  \]
  and for $r \in (0,1/2)$
  \[
    1+\frac{r}{2} \le \frac{M}{r} \le 1+r.
  \]
\end{remark}

\begin{figure}[!ht]
  \begin{center}
    \includegraphics[height=50mm]{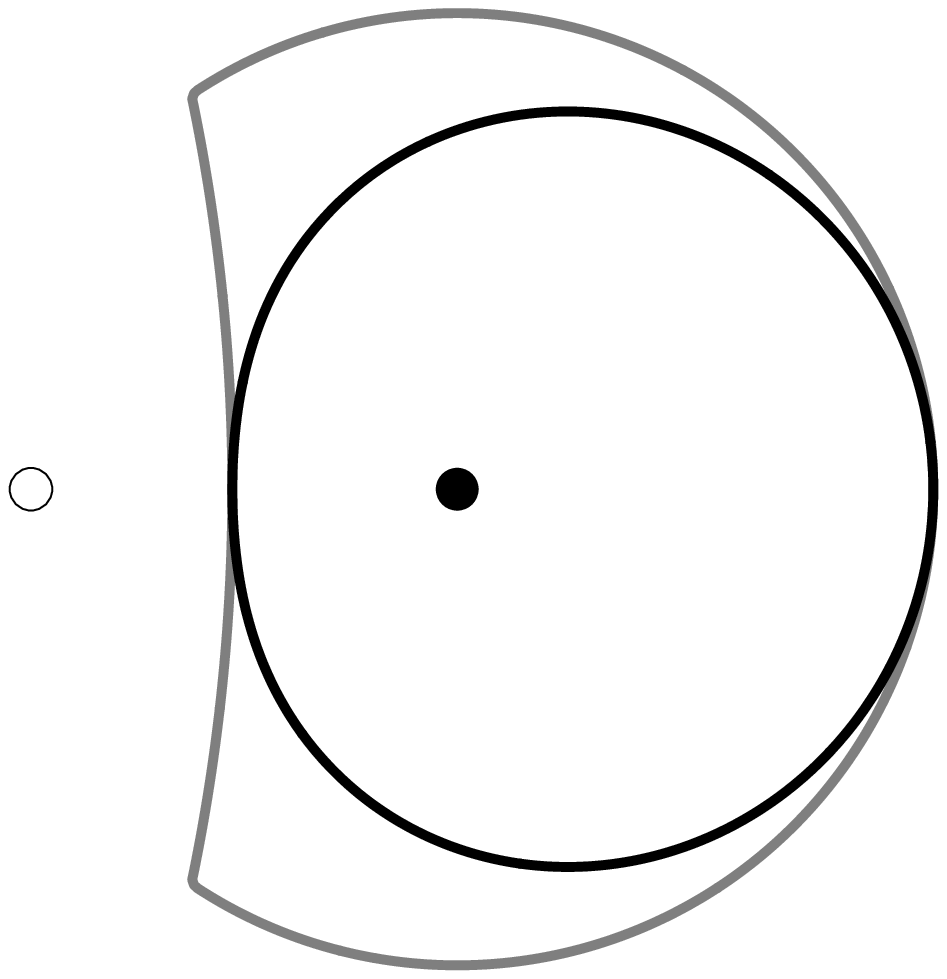}\hspace{1cm}
    \includegraphics[height=50mm]{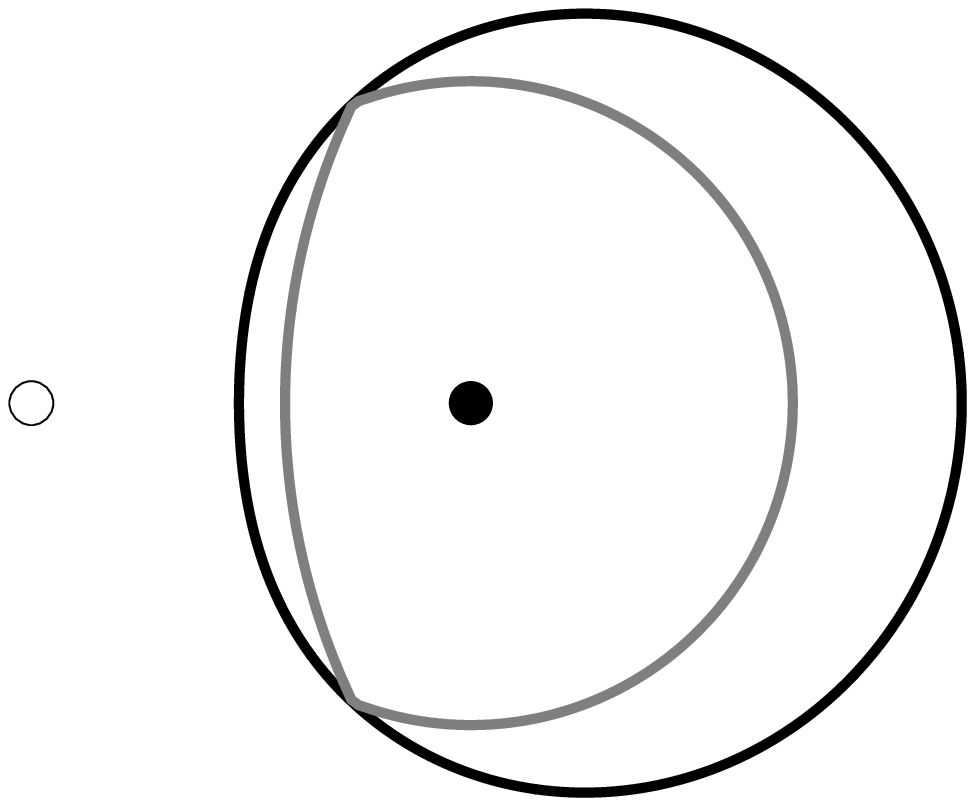}
    \caption{An example of inclusions of quasihyperbolic disks (black) and $j$-metric disks (gray) in punctured plane. The black dot is the center of the disks and the small circle denotes the origin. The radius of the quasihyperbolic disks is 0.75. \label{fig1}}
  \end{center}
\end{figure}

Results of Theorem \ref{jkj} and Corollary \ref{kjk} are illustrated in the planar case in Figure \ref{fig1}.

\begin{remark}
  We show that the radius $m_1 = \log ( 1+2 \sin (r/2) )$ of Theorem \ref{jkj} is better than the radius $m_2 = \log (2-e^{-r})$ of Proposition \ref{jkjkgen}, namely we show that the function
  \[
    f(r) = \log ( 1+2 \sin (r/2) ) - \log (2-e^{-r}) = \log \frac{1+2 \sin (r/2)}{2-e^{-r}}
  \]
  is positive on $(0,\log 3)$. The inequality $f(r) > 0$ is equivalent to $(1+2 \sin (r/2))/(2-e^{-r}) > 1$ and therefore it is sufficient to show that $g(r) = (1+2 \sin (r/2))/(2-e^{-r})$ is increasing and $g(0) = 1$. Since $g'(r) > 0$ is equivalent to $h(r)>0$ for $h(r) = (2e^r-1)\cos (r/2) - 2(1+ \sin(r/2))$ and $h(0)=0$, we need to show that $h'(r) > 0$. Because $h'(r) > 0$ is equivalent to $2 > \tan (r/2)$, it is true for $r \in (0,\log 3)$. Clearly $g(0) = 1$ and the assertion follows.
\end{remark}

We will next consider inclusion of the $j$-metric and the chordal metric balls.

\begin{theorem}\label{jqj}
  Let $G = \Rn \setminus \{ 0 \}$, $x \in G$ and $r \in (0,r_q(x))$. Then
  \[
    B_j(x,m) \subset B_q(x,r) \subset B_j(x,M),
  \]
  where
  \[
    m = \log \left( 1+r \left( |x|+\frac{1}{|x|} \right) \right)
  \]
  and
  \[
    M = \left\{ \begin{array}{ll} \displaystyle \log \frac{|x|(1-r^2(1+|x|^2))}{|x|-r\sqrt{1-r^2}(1+|x|^2)}, & \textrm{ for }|x| \le 1\\ \displaystyle \log \frac{|x|+r \sqrt{1-r^2}(1+|x|^2)}{|x|(1-r^2(1+|x|^2))}, & \textrm{ for }|x|>1. \end{array} \right.
  \]
  Moreover, the inclusions are sharp and $M/m \to 1$ as $r \to 0$.
\end{theorem}
\begin{proof}
  Let us first show that $B_j(x,m) \subset B_q(x,r)$. We assume $y \in B_j(x,m)$, which is equivalent to
  \[
    \frac{|x-y|}{\min \{ |x|,|y| \}} < r\left( |x|+\frac{1}{|x|}\right),
  \]
  and thus
  \[
    q(x,y)^2 < \frac{r^2 \min \{ |x|^2,|y|^2 \} \left( |x|+\frac{1}{|x|}\right)^2}{(1+|x|^2)(1+|y|^2)} = r^2 \min \left\{ 1,\frac{|y|^2}{|x|^2} \right\} \frac{(1+|x|^2)}{(1+|y|^2)} \le r^2.
  \]
  Therefore $y \in B_q(x,r)$.

  Let us then show that $B_q(x,r) \subset B_j(x,M)$. Since $r < 1/\sqrt{1+|x|^2}$, we have by Lemma \ref{qball}
  \begin{equation}\label{qballinjball1}
    B_q(x,r) = B^n(y,s),
  \end{equation}
  where $y = c x$, $c=c(x,r)  > 1$, and $s$ is depending only on $x$ and $r$. By (\ref{qballinjball1}) and \cite[p. 285]{k1} it is clear that for minimal $M$ such that $B_q(x,r) \subset B_j(x,M)$ we have $\partial B_q(x,r) \cap \partial B_j(x,M) \subset \{ z \in G \colon z = x t, t > 0 \}$. In other words it is sufficient to show that
  \begin{equation}\label{jballestimate}
    |y|+s \le |x|e^M \quad \textrm{and} \quad |y|-s \ge |x|e^{-M}.
  \end{equation}
  Before proving \eqref{jballestimate} we observe that since
  \[
    |y|^2-s^2  = \frac{|x|^2-r^2-r^2|x|^2}{1-r^2-r^2|x|^2}
  \]
  we have that $|x| \le 1$ is equivalent to $|x|/(|y|-s) \ge (|y|+s)/|x|$ and $|x| > 1$ is equivalent to $|x|/(|y|-s) < (|y|+s)/|x|$. Now we have for $|x| \le 1$
  \[
    \frac{|y|+s}{|x|} \le \frac{|x|}{|y|-s} = \frac{|x|(1-r^2(1+|x|^2))}{|x|-r\sqrt{1-r^2}(1+|x|^2)} = e^M
  \]
  and for $|x| > 1$
  \[
    \frac{|x|}{|y|-s} < \frac{|y|+s}{|x|} = \frac{|x|+r \sqrt{1-r^2}(1+|x|^2)}{|x|(1-r^2(1+|x|^2))} = e^M
  \]
  and the inequalities of (\ref{jballestimate}) hold true. 

  We show that $m$ is sharp. We choose $y \in G$ with $|y| = |x|$ and $q(x,y) = r$. Now $|x-y| = r(1+|x|^2)$ and
  \[
    j(x,y) = \log \left( 1+\frac{|x-y|}{|x|} \right) = \log \left( 1+\frac{r(1+|x|^2)}{|x|} \right) = m.
  \]

  We show that $M$ is sharp. Let us first assume $|x| \le 1$. We choose $y \in G$ with $q(x,y) = r$ and $|y| \le |z|$ for all $z \in G$ with $q(x,z) = r$. Now
  \[
    q(x,y) = \frac{|x|-|y|}{\sqrt{1+|x|^2}\sqrt{1+|y|^2}} = r
  \]
  implying
  \begin{equation}\label{mody}
    |y| = \frac{|x|-r\sqrt{1-r^2}(1+|x|^2)}{1-r^2(1+|x|^2)}.
  \end{equation}
  Now $j(x,y) = \log (|x|/|y|) = M$ and $M$ is sharp.

  Let us then assume $|x| > 1$. We choose $y \in G$ with $q(x,y) = r$ and $|y| \ge |z|$ for all $z \in G$ with $q(x,z) = r$ by (\ref{mody}) we have $j(x,y) = \log (|y|/|x|) = M$.

  Finally, we show that $M/m \to 1$ as $r \to 0$. Let us first consider the case $|x| \le 1$. By l'H\^opital's rule
  \[
    \lim_{r \to 0} \frac{M}{m} = \lim_{r \to 0} \frac{r+|x|+r|x|^2}{\sqrt{1-r^2}(1-2r^2)|x|+r(1-r^2)(1-|x|^2)} = 1.
  \]

  In the case $|x| > 1$ we obtain by l'H\^opital's rule
  \[
    \lim_{r \to 0} \frac{M}{m} = \lim_{r \to 0} \frac{r+|x|+r|x|^2}{\sqrt{1-r^2}(1-2r^2)|x|-r(1-r^2)(1-|x|^2)} = 1.
  \]
\end{proof}

\begin{corollary}\label{qjq}
  Let $G = \Rn \setminus \{ 0 \}$, $x \in G$ and $r \in (0,\log(1+|x|+1/|x|))$. Then
  \[
    B_q(x,m) \subset B_j(x,r) \subset B_q(x,M),
  \]
  where
  \[
    m = \left\{ \begin{array}{ll} \displaystyle \frac{(e^r-1)|x|}{\sqrt{(1+|x|^2)(e^{2r}+|x|^2)}}, & \textrm{if }|x| \le 1,\\ \displaystyle \frac{(e^r-1)|x|}{\sqrt{(1+|x|^2)(1+|x|^2 e^{2r})}}, & \textrm{if } |x| \ge 1, \end{array} \right.
  \]
  and
  \[
    M = \frac{|x|(e^r-1)}{1+|x|^2}.
  \]
  Moreover, the inclusions are sharp and $M/m \to 1$ as $r \to 0$. The assumption $r < \log(1+|x|+1/|x|)$ is equivalent to $M<1$.
\end{corollary}
\begin{proof}
  The claim follows from Theorem \ref{jqj}.
\end{proof}

\begin{figure}[!ht]
  \begin{center}
    \includegraphics[height=50mm]{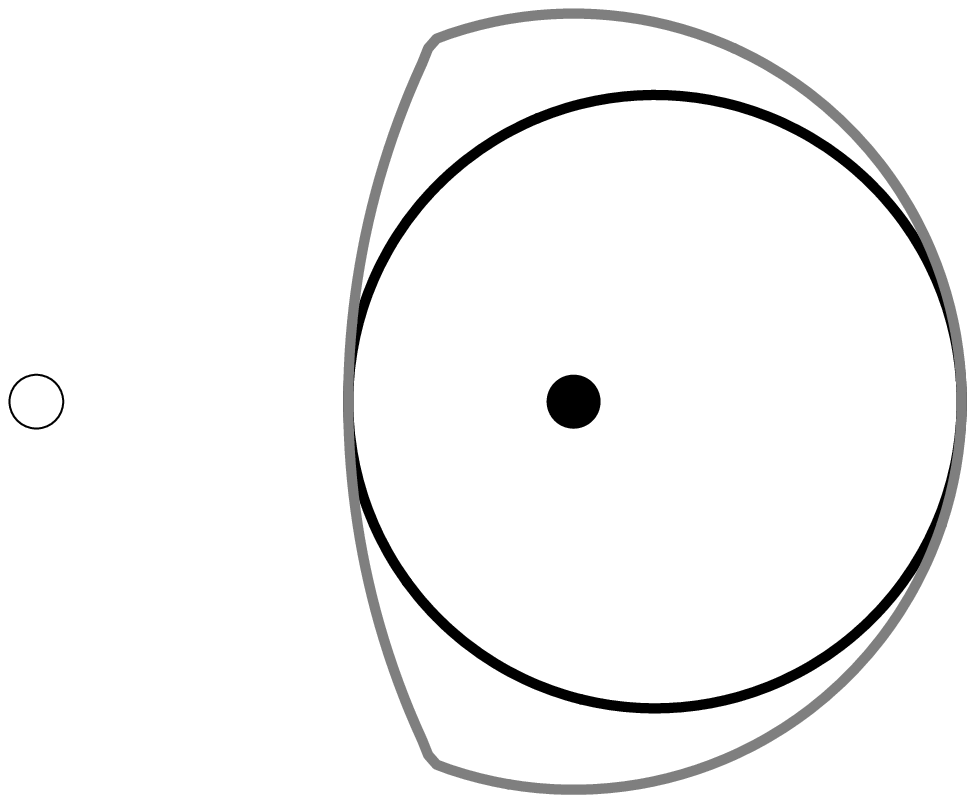}\hspace{1cm}
    \includegraphics[height=50mm]{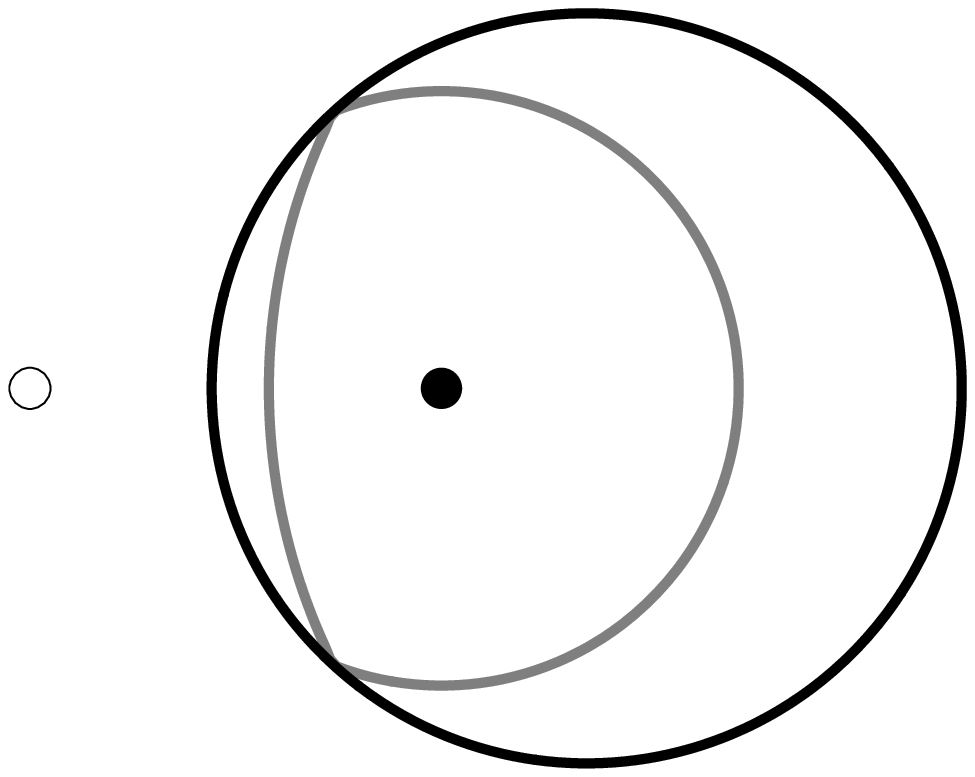}
    \caption{An example of inclusions of chordal disks (black) and $j$ metric disks (gray) in punctured plane. The black dot is the center of the disks and the small circle denotes the origin. The radius of the $j$-metric disks is 0.5. \label{fig2}}
  \end{center}
\end{figure}

The results of Theorem \ref{jqj} and Corollary \ref{qjq} are illustrated in the planar case in Figure \ref{fig2}.

We will finally consider inclusion of the quasihyperbolic and the chordal metric balls.

\begin{theorem}\label{kqk}
  Let $G = \Rn \setminus \{ 0 \}$, $x \in G$ and $r \in (0,R)$ for
  \[
    R =  \frac{2|x|}{\sqrt{1+|x|^2}\sqrt{1+9|x|^2}}.
  \]
  Then
  \[
    B_k(x,m) \subset B_q(x,r) \subset B_k(x,M),
  \]
  where
  \[
    m = \log \left( 1+r \left( |x|+\frac{1}{|x|} \right) \right)
  \]
  and
  \[
    M = 2 \arcsin f(r,x)
  \]
  with
  \[
    f(r,x) = \left\{ \begin{array}{ll} \displaystyle \frac{r(1+|x|^2)}{2\sqrt{1-r^2}|x|-2r}, & \textrm{ for }|x| \le 1\\ \displaystyle \frac{r+r|x|^2}{2\sqrt{1-r^2}|x|-2r|x|^2}, & \textrm{ for }|x|>1. \end{array} \right.
  \]
  Moreover, we have $M/m \to 1$ as $r \to 0$.
\end{theorem}
\begin{proof}
  By \eqref{gehringpalka} and Theorem \ref{jqj} we have
  \[
    B_k(x,m) \subset B_j(x,m) \subset B_q(x,r)
  \]
  for $r \in (0,r_q(x))$.
  
  By Theorem \ref{jqj} we have
  \begin{equation}\label{qj}
    B_q(x,r) \subset B_j(x,t)
  \end{equation}
  for $r \in (0,r_q(x))$. By Theorem \ref{jkj}
  \begin{equation}\label{jk}
    B_j(x,t) \subset B_k(x,s)
  \end{equation}
  for $t \in (0,\log 3)$ and $s = 2 \arcsin ((e^t-1)/2)$. Combining (\ref{qj}) and (\ref{jk}) gives $M = s$ and $r \in (0,R_0) \subset (0,r_q(x))$, where
  \[
    R_0 =  \min \left\{ \frac{|x|(e^{\pi/2}-1)}{\sqrt{1+|x|^2}\sqrt{1+e^\pi |x|^2}} , \frac{2|x|}{\sqrt{1+|x|^2}\sqrt{1+9|x|^2}} \right\}.
  \]
  We also have $f(r,x) \in [0,1]$ for $r \in (0,R_0)$.

  We show next that $R=R_0$, which is equivalent to
  \begin{equation}\label{RandR0}
    \frac{e^{\pi/2}-1}{\sqrt{1+e^\pi |x|^2}} \ge \frac{2}{\sqrt{1+9 |x|^2}}.
  \end{equation}
  Inequality \eqref{RandR0} is equivalent to $|x|^2 (5e^{\pi/2}-3) \ge -(e^{\pi/2}+1)$, which is clearly true and thus $R=R_0$.

  Finally, by a straightforward computation using l'H\^opital's  rule we obtain $M/m \to 1$ as $r \to 0$.
\end{proof}

\begin{corollary}\label{qkq}
  Let $G = \Rn \setminus \{ 0 \}$, $x \in G$ and $r \in (0,\log(1+|x|+1/|x|))$. Then
  \[
    B_q(x,m) \subset B_k(x,r) \subset B_q(x,M),
  \]
  where
  \begin{eqnarray*}
    m & = & \min \left\{ \frac{2|x|\sin (r/2)}{\sqrt{(1+|x|^2)(|x|^2+(1+2\sin(r/2))^2)}} \right.,\\ && \left. \frac{2|x|\sin (r/2)}{\sqrt{(1+|x|^2)(1+|x|^2(1+2\sin(r/2))^2))}} \right\}
  \end{eqnarray*}
  and
  \[
    M = \frac{|x|(e^r-1)}{1+|x|^2}.
  \]
  Moreover, we have $M/r \to 1$ as $r \to 0$.
\end{corollary}
\begin{proof}
  The assertion follows from Theorem \ref{kqk}. The assumption $r < \log(1+|x|+1/|x|)$ is equivalent to $M<1$.
\end{proof}

Note that the radii of Theorem \ref{kqk} and Corollary \ref{qkq} are not sharp. Radii $m$ and $M$ of Theorem \ref{jqj}, Corollary \ref{qjq}, Theorem \ref{kqk} and Corollary \ref{qkq} depend on $r$ and $|x|$. We can use \eqref{boundsforx} to make some of the radii independent of $|x|$.

\begin{theorem}\label{main1}
  Let $G = \Rn \setminus \{ 0 \}$, $x \in G$ and $r \in (0,r_q(x))$. Then
  \[
    B_j(x,m) \subset B_q(x,r) \quad \textrm{and} \quad B_k(x,m) \subset B_q(x,r),
  \]
  where
  \[
    m = \log \left( 1+\frac{2r^2}{\sqrt{1-r^2}} \right).
  \]
\end{theorem}
\begin{corollary}
  Let $G = \Rn \setminus \{ 0 \}$, $x \in G$ and $r \in (0,\log(1+\sqrt{2}))$. Then
  \[
    B_j(x,r) \subset B_q(x,M) \quad \textrm{and} \quad B_k(x,r) \subset B_q(x,M),
  \]
  where
  \[
    M = \frac{\sqrt{(e^r-1)(e^r-1+\sqrt{17+e^r(e^r-2)})}}{2\sqrt{2}}.
  \]
\end{corollary}
Note that for $M$ such that $B_q(x,r) \subset B_j(x,M)$ we have by Theorem \ref{jqj} that $M \to \infty$ whenever $|x| \to r/\sqrt{1-r^2}$ or $|x| \to \sqrt{1-r^2}/r$. Similarly by Theorem \ref{kqk} we have that $M \to \infty$ whenever $|x| \to r/\sqrt{1-r^2}$ or $|x| \to \sqrt{1-r^2}/r$ for $M$ such that $B_q(x,r) \subset B_k(x,M)$.

\begin{remark}
  In this paper we assume that metric balls have the same center point. Similar problems can be considered without this assumption. For example it can be proved that for $x \in G = \PS$, $r > 0$
  \[
    \sup \{ t \colon B^n(z,t) \subset B_k(x,r),\, z \in G \} = \cosh r.
  \]
\end{remark}

\section{Half-space}\label{halfspace}

In this section we will consider inclusion of metric balls defined by the metrics $q$, $j_G$ and $k_G$ for $G = \Hn = \{  x \in \Rn \colon x_n > 0 \}$. Recall that in $\Hn$ the quasihyperbolic metric agrees with the hyperbolic metric
\begin{equation}\label{hypmetric}
    \cosh k_\Hn(x,y) = \cosh \rho_\Hn(x,y) = 1+\frac{|x-y|^2}{2x_ny_n} \quad \textrm{for } x,y \in \Hn.
\end{equation}

\begin{theorem}\label{jkj2}
  Let $G = \Hn$, $x \in G$ and $r > 0$. Then
  \[
    B_j(x,m) \subset B_k(x,r),
  \]
  where
  \[
    m = \log \left( 1+\sqrt{2} \sqrt{\cosh r-1} \right).
  \]
  Moreover, the inclusion is sharp and $r/m \to 1$ as $r \to 0$.
\end{theorem}
\begin{proof}
  Let $y \in B_j(x,m)$ and denote $y = y_n e_n + y'$. By \cite[(2.11)]{vu3} we have $B_k(x,r) = B^n(z,|x| \sinh r)$ for $z = |x| e_n \cosh r$. Let us first assume $x_n = d(x) \le d(y) = y_n$. Now $y \in B_j(x,m)$ implies that $|y'|^2 < 2(\cosh r -1)x_n^2-(x_n-y_n)^2$. Now
  \begin{eqnarray*}
    |z-y|^2 & = & (x_n \cosh r-y_n)^2+|y'|^2\\
    & < & (x_n \cosh r-y_n)^2+2(\cosh r -1)x_n^2-(x_n-y_n)^2\\
    & = & x_n^2(\cosh^2r+2\cosh r-3)-2x_n y_n (\cosh r-1)\\
    & \le & x_n^2(\cosh^2r+2\cosh r-3)-2x_n^2 (\cosh r-1)\\
    & = & x_n^2(\cosh^2r-1)\\
    & \le & x_n^2 \sinh^2 r
  \end{eqnarray*}
  and therefore $y \in B_k(x,r)$.

  Let us then assume $y_n = d(y) \le d(x) = x_n$. Now $y \in B_j(x,m)$ implies that $|y'|^2 < 2(\cosh r -1)y_n^2-(x_n-y_n)^2$. Now
  \begin{eqnarray*}
    |z-y|^2 & = & (x_n \cosh r-y_n)^2+|y'|^2\\
    & < & (x_n \cosh r-y_n)^2+2(\cosh r -1)y_n^2-(x_n-y_n)^2\\
    & = & x_n^2(\cosh^2r-1)+2(\cosh r-1)y_n(y_n-x_n)\\
    & \le & x_n^2(\cosh^2r-1)\\
    & = & x_n^2 \sinh^2 r
  \end{eqnarray*}
  and therefore $y \in B_k(x,r)$.

  We show that $m$ is sharp for $y \in \Hn$ such that $j(x,y) = m$ and $d(x)= d(y)$. Now $|y'|^2 = 2(\cosh r -1)x_n^2$ and
  \[
    |z-y|^2 = x_n^2(\cosh r-1)^2+|y'|^2 = x_n^2(\cosh r-1)^2+2(\cosh r -1)x_n^2 = x_n^2 \sinh^2 r.
  \]

  Finally, by l'H\^opital's  rule
  \begin{eqnarray*}
    \lim_{r \to 0} \frac{r}{m} & = & \lim_{r \to 0} \frac{2\cosh r-2-\sqrt{2}\sqrt{\cosh r-1}}{\sinh r}\\
    & = & \lim_{r \to 0} \frac{\displaystyle 2\sinh r+\frac{\sinh r}{\sqrt{2}\sqrt{\cosh r-1}}}{\cosh r}\\
    & = & \lim_{r \to 0} \frac{\sinh r}{\sqrt{2}\sqrt{\cosh r-1}} = \lim_{r \to 0} \frac{\sqrt{\cosh r+1}}{\sqrt{2}} = 1
  \end{eqnarray*}
  and the assertion follows.
\end{proof}

\begin{corollary}\label{kjk2}
  Let $G = \Hn$, $x \in G$ and $r > 0$. Then
  \[
    B_j(x,r) \subset B_k(x,M),
  \]
  where
  \[
    M = \arcosh \left( 1+\frac{(e^r -1)^2}{2} \right).
  \]
  Moreover, the inclusion is sharp and $M/r \to 1$ as $r \to 0$.
\end{corollary}
\begin{proof}
  The assertion follows from Theorem \ref{jkj2}.
\end{proof}

\begin{remark}
  Let $m$ as in Theorem \ref{jkj2}. It is possible to show that for $r \in (0,1)$
  \[
    1+\frac{2r}{5} \le \frac{r}{m} \le 1+\frac{r}{2} \quad \textrm{and} \quad 1+\frac{r}{2} \le \frac{M}{r} \le 1+\frac{3r}{5}.
  \]
\end{remark}

\begin{figure}[!ht]
  \begin{center}
    \includegraphics[height=40mm]{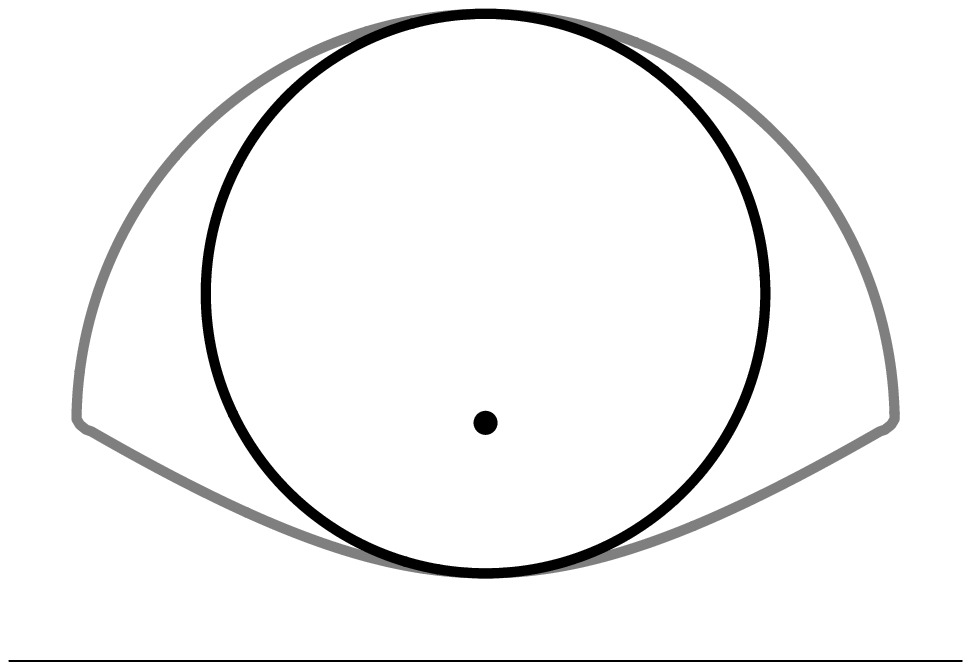}\hspace{5mm}
    \includegraphics[height=40mm]{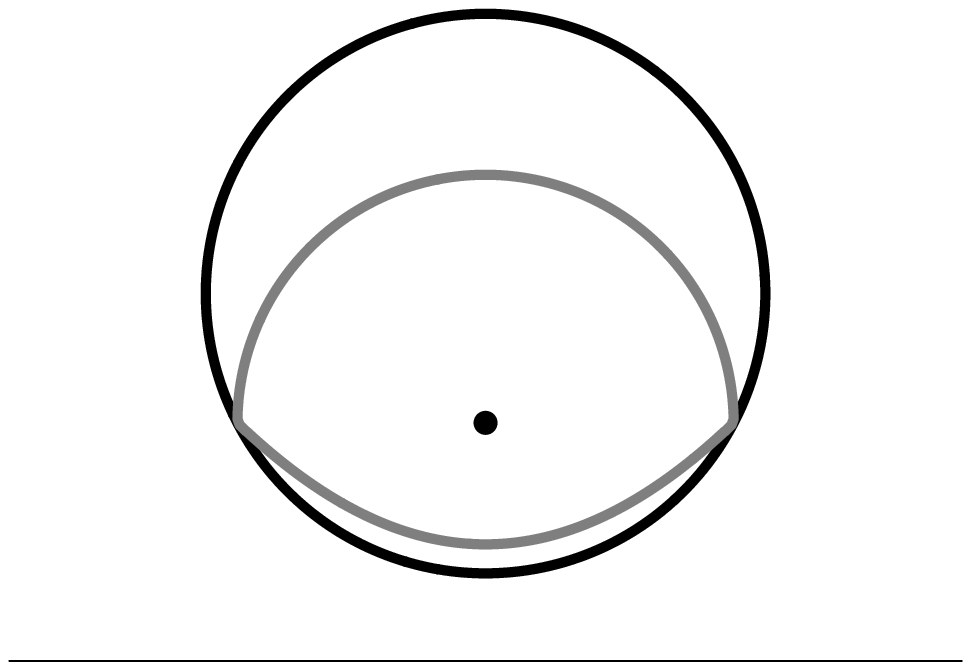}
    \caption{An example of inclusions of quasihyperbolic disks (black) and $j$ metric disks (gray) in half-plane. The black line represents $\partial G$. The radius of the quasihyperbolic disks is 1.\label{fig3}}
  \end{center}
\end{figure}

Results of Theorem \ref{jkj2} and Corollary \ref{kjk2} are illustrated in the planar case in Figure \ref{fig3}.

\begin{remark}
  1) We show that the radius $m_1 = \log ( 1+\sqrt{2} \sqrt{\cosh r-1} )$ of Theorem \ref{jkj2} is better than the radius $m_2 = \log (2-e^{-r})$ of Proposition \ref{jkjkgen}, namely we show that
  \begin{equation}\label{remarkcomp1}
    1+ \displaystyle \sqrt{2} \sqrt{\cosh r-1} > 2-e^{-r}
  \end{equation}
  for $r \in (0,\infty)$. Equation \eqref{remarkcomp1} is equivalent to $e^{2r}+3 - (3e^r+e^{-r}) > 0$ and because $e^{2r}+3 - (3e^r+e^{-r}) = (e^r-1)^3 e^{-r}$ the claim is clear.

  2) We show that the radius $m_1 = \log ( 1+\sqrt{2} \sqrt{\cosh r-1} )$ of Theorem \ref{jkj2} is better than the radius $m_3 = \log (1+2 \sin (r/2))$ of Theorem \ref{jkj}. We show that
\begin{equation}\label{remarkcomp2}
    1+ \displaystyle \sqrt{2} \sqrt{\cosh r-1} \ge 1+2 \sin (r/2)
\end{equation}
  for $r \in (0,\infty)$. Inequality \eqref{remarkcomp2} is equivalent to $e^r+e^{-r} \ge 4\sin^2 (r/2)$, which true because $e^r+e^{-r} > 2+r^2$ and $4 \sin^2 (r/2) < r^2$.
\end{remark}

To simplify notation we may assume that $x \in \Hn$ with $x_1 = x_2 = \cdots = x_{n-1} = 0$. Since we also want that the chordal balls $B_q(x,r)$ are in $\Hn$ it is natural to assume that $r \in (0,r_q(x))$. This assumption is same as in the case $G = \PS$ and thus \eqref{boundsforx} holds.

\begin{theorem}\label{jqj2}
  Let $G = \Hn$, $x \in G$ with $x_1 = x_2 = \cdots = x_{n-1} = 0$ and $r \in (0,r_q(x))$. Then
  \[
    B_j(x,m) \subset B_q(x,r) \subset B_j(x,M),
  \]
  where
  \[
    m = \min \left\{ \log \left( 1+\frac{r (1+ |x|^2)}{|x|\sqrt{1-r^2-r^2 |x|^2}} \right) , \log \left( 1+\frac{r(1+|x|^2)}{\sqrt{1-r^2}|x|-r} \right) \right\}
  \]
  and
  \[
    M = \max \left\{ \log \left( 1+\frac{r(1+|x|^2)}{\sqrt{1-r^2}|x|-r} \right), \log \left( 1+\frac{r(1+|x|^2)}{|x|(\sqrt{1-r^2}-r|x|)} \right) \right\}.
  \]
  Moreover, the inclusions are sharp and $M/m \to 1$ as $r \to 0$.
\end{theorem}
\begin{proof}
  By Lemma \ref{qball} we know that $B_q(x,r) = B^n(y,s)$ for given $s$ and $y$.

  We start by proving the bound $m$. Let $z \in B_j(x,m)$. We assume first that $|x| = d(x) \le d(z)$ and show that in this case that $m = \log ( 1+(r (1+ |x|^2))/(|x|\sqrt{1-r^2-r^2 |x|^2})$. By definition
  \[
    |x-z| < (e^m-1)|x| = \frac{r(1+|x|^2)}{\sqrt{1-r^2-r^2|x|^2}}
  \]
  and therefore $z \in B^n(x,(e^m-1)|x|)$ and
  \[
    q(x,z) < q(x,z_1)
  \]
  for $z_1 = |x|e_n + |x|(e^m-1) e_1$. Now
  \begin{equation}\label{upperboundforq}
    q(x,z) < q(x,z_1) = \frac{|x-z_1|}{\sqrt{1+|x|^2}\sqrt{1+|z_1|^2}} = r
  \end{equation}
  and the assertion follows. This also shows that $m$ is sharp for $d(x) \le d(z)$.

  Let us then assume $d(z) \le d(x)$. Now for $z \in \partial B_j(x,m)$ we have $|x-z| = (e^m-1)d(z)$ and thus for each plane $p$, which contains $\{ z\in \Rn \colon z=te_n, t \in \R \}$, the intersection $\partial B_j(x,m) \cap p$ is a plane curve with largest curvature at the point that is closest to $\partial G$. Therefore
  \[
    q(x,z) \le \min \{ q(x,z_1),q(x,z_2) \}
  \]
  for $z_1 = |x|e_n + |x|(e^m-1) e_1$ and $z_2 = (|x|\sqrt{1-r^2}-r)/(\sqrt{1+r^2}+r|x|) e_n$. By \eqref{upperboundforq} we have $q(x,z) < q(x,z_1) = r$. Because
  \begin{equation}\label{zys}
    |z_2| = \frac{|x|\sqrt{1-r^2}-r}{\sqrt{1+r^2}+r|x|} \le \frac{|x|-r(1+|x|^2)\sqrt{1-r^2}}{1-r^2(1+|x|^2)} = |y|-s
  \end{equation}
  we have $q(x,z) < q(x,z_2) = r$ and the assertion follows. The second inequality in \eqref{zys} is equivalent to
  \[
    (1+|x|^2) (r^2 (2|x|\sqrt{1-r^2}-r)+r(1-\sqrt{1-r^4}))+|x|(\sqrt{1+r^2}-\sqrt{1-r^2}) > 0
  \]
  and holds true, because by assumption $2|x|\sqrt{1-r^2} > r$, $1 > \sqrt{1-r^4}$ and $\sqrt{1+r^2} > \sqrt{1-r^2}$. Sharpness of $m$ follows from the selection of $z_1$ and $z_2$.

  Let us then prove the bound $M$. We assume that $z \in \partial B_q(x,r)$ and $|x|=d(x)\le d(z)$. By definition $|x-z| = (e^M-1)|x|$ and therefore
  \[
    j(x,z) \le j(x,z_3) = \log \left( 1+\frac{|y|+s-|x|}{|x|} \right) = \log \left( 1+\frac{r(1+|x|^2)}{|x|(\sqrt{1-r^2}-r|x|)} \right)
  \]
  for $z_3 = (|y|+s) e_n$.

  Let us assume $z \in \partial B_q(x,r)$ and $d(z) \le d(x)=|x|$. By definition $|x-z| = (e^M-1)|z|$ and therefore
  \[
    j(x,z) \le j(x,z_2) = \log \left( 1+\frac{|y|-s-|x|}{|x|} \right) = \log \left( 1+\frac{r(1+|x|^2)}{\sqrt{1-r^2}|x|-r} \right)
  \]
  for $z_2 = (|y|-s) e_n$.

  Sharpness of $M$ is clear by selection of $z_3$ and $z_2$. We shall finally show that $M/m \to 1$ as $r\to 0$. We denote
  \[
    m_1 = \log \left( 1+\frac{r (1+ |x|^2)}{|x|\sqrt{1-r^2-r^2 |x|^2}} \right), \quad m_2 =  \log \left( 1+\frac{r(1+|x|^2)}{\sqrt{1-r^2}|x|-r} \right)
  \]
  and
  \[
    M_1 = \log \left( 1+\frac{r(1+|x|^2)}{\sqrt{1-r^2}|x|-r} \right), \quad M_2 = \log \left( 1+\frac{r(1+|x|^2)}{|x|(\sqrt{1-r^2}-r|x|)} \right).
  \]
  Clearly $M_1/m_2 = 1$ and by l'H\^opital's  rule
  \begin{eqnarray*}
    \lim_{r \to 0} \frac{M_1}{m_1} & = & \lim_{r \to 0} \frac{(r^2(1+|x|^2)-1)(r+r|x|^2+|x|\sqrt{1-r^2(1+|x|^2)})}{|x|\sqrt{1-r^2}(2r^2-1)+r(r^2-1)(|x|^2-1)} = 1,\\
    \lim_{r \to 0} \frac{M_2}{m_1} & = & \lim_{r \to 0} \frac{-(r^2(1+|x|^2)-1)(r+r|x|^2+|x|\sqrt{1-r^2(1+|x|^2)})}{-|x|\sqrt{1-r^2}(2r^2-1)+r(r^2-1)(|x|^2-1)} = 1,\\
    \lim_{r \to 0} \frac{M_2}{m_2} & = & \lim_{r \to 0} \frac{|x|-2r^2|x|+r\sqrt{1-r^2}(|x|^2-1)}{|x|-2r^2|x|-r\sqrt{1-r^2}(|x|^2-1)} = 1
  \end{eqnarray*}
  and the assertion follows.
\end{proof}

\begin{corollary}\label{qjq2}
  Let $G = \Hn$, $x \in G$ with $x_1 = x_2 = \cdots = x_{n-1} = 0$ and $r > 0$. Then
  \[
    B_q(x,m) \subset B_j(x,r) \subset B_q(x,M),
  \]
  where
  \[
    m =\left\{ \begin{array}{ll} \displaystyle \frac{|x|(e^r-1)}{\sqrt{1+|x|^2}\sqrt{e^{2r}+|x|^2}}, & \textrm{if }|x| \le 1,\\ \displaystyle \frac{|x|(e^r-1)}{\sqrt{1+|x|^2}\sqrt{1+e^{2r}|x|^2}}, & \textrm{if }|x| \ge 1 \end{array} \right.
  \]
  and
  \[
    M = \left\{ \begin{array}{ll} \displaystyle \frac{|x|(e^r-1)}{\sqrt{1+|x|^2}\sqrt{e^{2r}+|x|^2}}, & \textrm{if } |x| \le \sqrt{\tanh (r/2)}, \\ \displaystyle \frac{|x|(e^r-1)}{\sqrt{1+|x|^2}\sqrt{1+(2+e^r(e^r-2))|x|^2}}, & \textrm{if }|x| \ge \sqrt{\tanh (r/2)}. \end{array} \right.
  \]
  Moreover, the inclusions are sharp and $M/m \to 1$ as $r \to 0$.
\end{corollary}
\begin{proof}
  The assertion follows from Theorem \ref{jqj2}.
\end{proof}

\begin{figure}[!ht]
  \begin{center}
    \includegraphics[height=45mm]{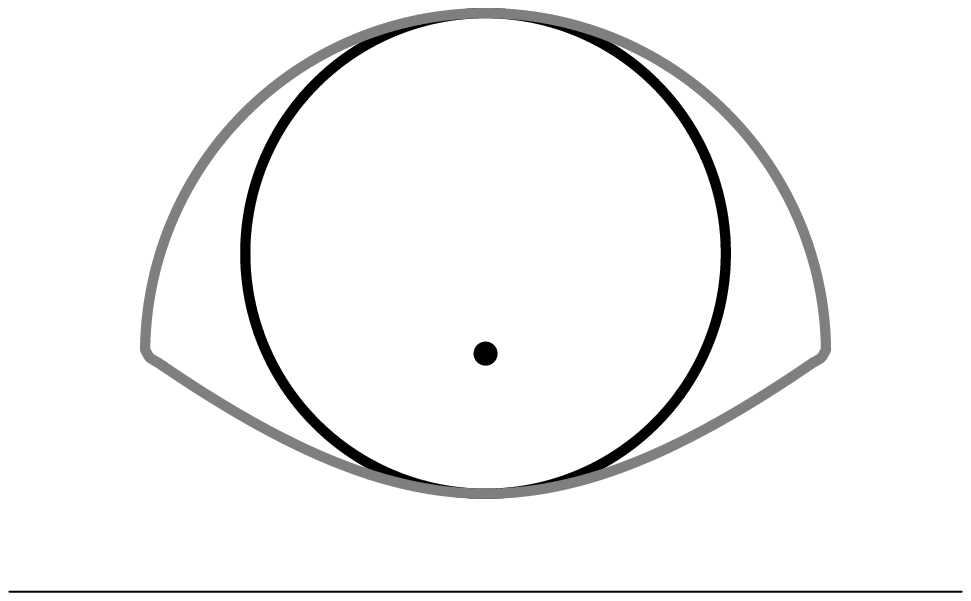}\hspace{5mm}
    \includegraphics[height=45mm]{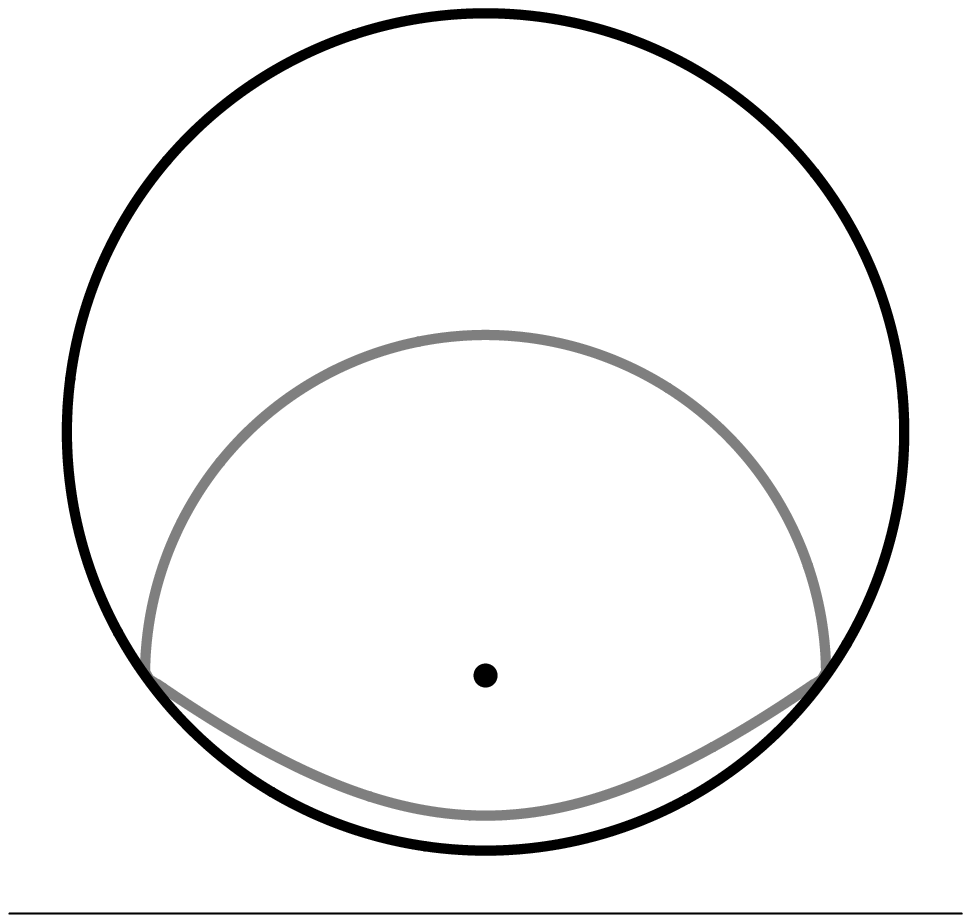}
    \caption{An example of inclusions of chordal disks (black) and $j$ metric disks (gray) in half-plane. The black line represents $\partial G$. The radius of the $j$-metric disks is 0.9. \label{fig4}}
  \end{center}
\end{figure}

Results of Theorem \ref{jqj2} and Corollary \ref{qjq2} are illustrated in the planar case in Figure \ref{fig4}.

\begin{remark}
  We show that Theorem \ref{jqj2} improves Theorem \ref{jqj}. Firstly, we consider the radius $m$ and show that
  \[
    f(r) = \min \left\{ \frac{1}{\sqrt{1-r^2-r^2|x|^2}} , \frac{|x|}{\sqrt{1-r^2}|x|-r} \right\} - 1
  \]
  is positive for $r \in (0,r_q(x))$. If $|x| \le 1$, then $f(r) = 1/\sqrt{1-r^2-r^2|x|^2}$ and $f(r) > 0$ is equivalent to $r^2(1+|x|^2) > 0$, which is true. If $|x| > 1$, then $f(r) = |x|/(\sqrt{1-r^2}|x|-r)$ and $f(r) > 0$ is equivalent to $|x| > -r/(1-\sqrt{1-r^2})$, which is true because $-r/(1-\sqrt{1-r^2}) < 0$. 

  Secondly, we consider the radius $M$. It is possible to show that the value $M$ in Theorem \ref{jqj2} is less than or equal to the value $M$ in Theorem \ref{jqj}. The calculation is straightforward but tedious and we omit it.
\end{remark}

\begin{theorem}\label{kqk2sharp}
  Let $G = \Hn$, $x \in G$ with $x_1 = x_2 = \cdots = x_{n-1} = 0$ and $r \in (0,r_q(x))$. Then
  \[
    B_k(x,m) \subset B_q(x,r) \subset B_k(x,M),
  \]
  where
  \[
    m = \min \left\{ \log \frac{|x|+r \sqrt{1-r^2}(1+|x|^2)}{|x|(1-r^2(1+|x|^2))} , \log \frac{|x|(r^2(1+|x|^2)-1)}{r \sqrt{1-r^2}(1+|x|^2)-|x|} \right\}
  \]
  and
  \[
    M = \max \left\{ \log \frac{|x|+r \sqrt{1-r^2}(1+|x|^2)}{|x|(1-r^2(1+|x|^2))} , \log \frac{|x|(r^2(1+|x|^2)-1)}{r \sqrt{1-r^2}(1+|x|^2)-|x|} \right\}.
  \]
  Moreover, the inclusions are sharp and $M/m \to 1$ as $r \to 0$.
\end{theorem}
\begin{proof}
  Let us denote
  \[
    a = \log \frac{|x|+r \sqrt{1-r^2}(1+|x|^2)}{|x|(1-r^2(1+|x|^2))} \quad \textrm{and} \quad A = \log \frac{|x|(r^2(1+|x|^2)-1)}{r \sqrt{1-r^2}(1+|x|^2)-|x|}.
  \]
  It is easy to verify that $a=A$ is equivalent to $|x| = 1$, $a<A$ is equivalent to $|x| < 1$ and $a>A$ is equivalent to $|x| > 1$.

  Since $k$ agrees with the usual hyperbolic metric in $G$ we know that quasihyperbolic balls are Euclidean balls. Therefore, both $B_k(x,r)$ and $B_q(x,r)$ are Euclidean balls of form $B^n(s x_n,t)$ for $s,t > 0$. For chordal metric the value of $s$ is $s_1 = |x|/(1-r^2(1+|x|^2))$ and quasihyperbolic metric $s_2 = |x| \cosh a$ and $s_3 = |x| \cosh A$. We will first show that for $|x| < 1$ we have $s_1 < s_2 < s_3$ and for $|x| > 1$ we have $s_3 < s_2 < s_1$. In the case $|x| = 1$ we have $s_1 = s_2 = s_3$ and the assertion follows.

Let us first assume that $|x| < 1$. Since $\cosh t$ is strictly increasing on $(0,\infty)$ and $a < A$ we have $s_2 < s_3$. Inequality $s_1 < s_2$ is equivalent to
\begin{equation}\label{fxr}
  2 f(x,r) < 1+f(x,r)^2,
\end{equation}
for
\[
  f(x,r) = \frac{|x|+r\sqrt{1-r^2}(1+|x|^2)}{|x|(1-r^2(1+|x|^2))}.
\]
By selection of $r$ we have $f(x,r) > 0$ and \eqref{fxr} implies that $s_1 < s_2$.

Let us now prove the inclusion $B_k(x,m) \subset B_q(x,r)$. Since $s_1 < s_2$ we have for the largest possible $m$ that $\partial B_k(x,m) \cap \partial B_q(x,r) = \{ y \}$, where $y \in \partial B_q(x,r)$ is such that $d(y) > d(w)$ for all $w \in \partial B_q(x,r)$, $w \neq y$. Now
\[
  q(x,y) = k(x,y)
\]
and thus $B_k(x,m) \subset B_q(x,r)$ and $m$ is sharp.

Let us then prove the inclusion $B_q(x,r) \subset B_k(x,M)$. Since $s_1 < s_3$ we have for the largest possible $m$ that $\partial B_k(x,m) \cap \partial B_q(x,r) = \{ z \}$, where $z \in \partial B_q(x,r)$ is such that $d(z) < d(w)$ for all $w \in \partial B_q(x,r)$, $w \neq z$. Now
\[
  q(x,z) = k(x,z)
\]
and thus $B_q(x,r) \subset B_k(x,M)$ and $M$ is sharp.

The case $|x| > 1$ is similar to the case $|x| < 1$.

We will finally show that $M/m \to 1$ as $r \to 0$. By l'H\^opital's  rule
\[
  \lim_{r \to 0} \frac{M}{m} = \lim_{r \to 0} \frac{|x|-2r^2|x|-r\sqrt{1-r^2}(|x|^2-1)}{|x|-2r^2|x|+r\sqrt{1-r^2}(|x|^2-1)} = 1.
\]
\end{proof}

\begin{corollary}\label{qkq2sharp}
  Let $G = \Hn$, $x \in G$ with $x_1 = x_2 = \cdots = x_{n-1} = 0$ and $r > 0$. Then
  \[
    B_q(x,m) \subset B_k(x,r) \subset B_q(x,M),
  \]
  where
  \[
    m = \min \left\{ \frac{|x|(e^r-1)}{\sqrt{(1+|x|^2)(e^{2r}+|x|^2)}} , \frac{|x|(e^r-1)}{\sqrt{(1+|x|^2)(1+e^{2r}|x|^2)}} \right\}
  \]
  for
  \[
    M = \max \left\{ \frac{|x|(e^r-1)}{\sqrt{(1+|x|^2)(e^{2r}+|x|^2)}} , \frac{|x|(e^r-1)}{\sqrt{(1+|x|^2)(1+e^{2r}|x|^2)}} \right\}.
  \]
  Moreover, the inclusions are sharp and $M/m \to 1$ as $r \to 0$.
\end{corollary}
\begin{proof}
  The assertion follows from Theorem \ref{kqk2sharp}.
\end{proof}

\begin{figure}[!ht]
  \begin{center}
    \includegraphics[height=50mm]{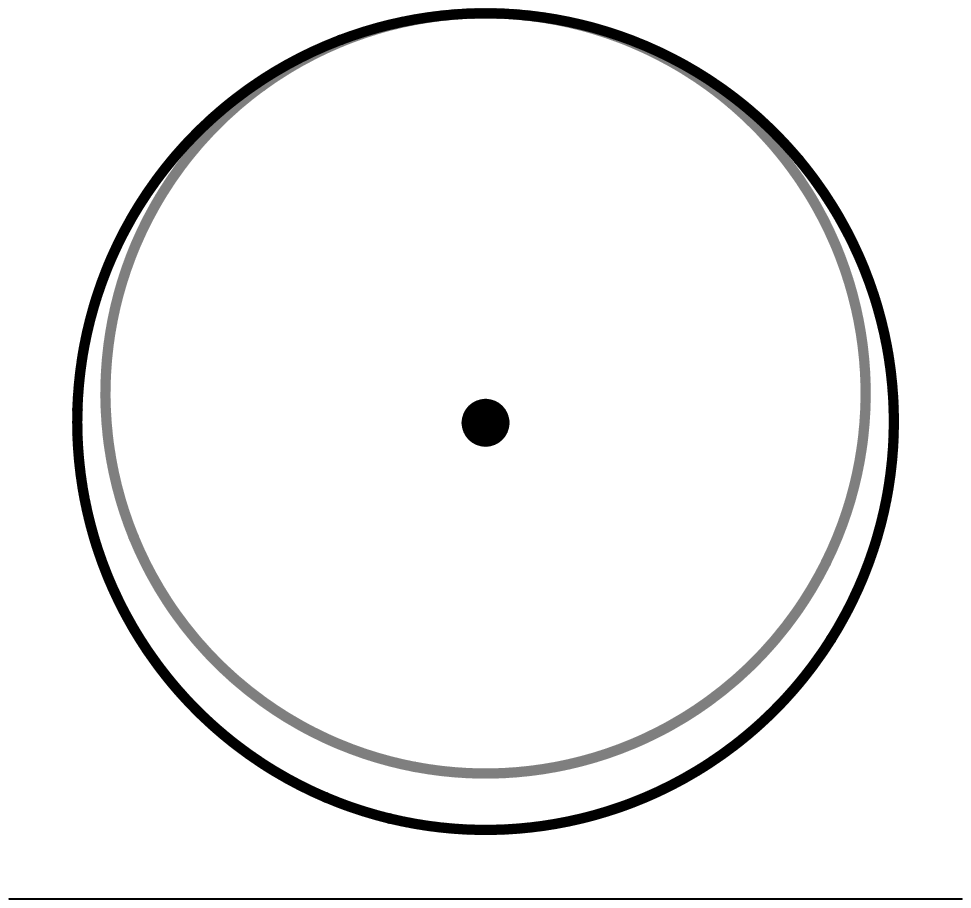}\hspace{5mm}
    \includegraphics[height=50mm]{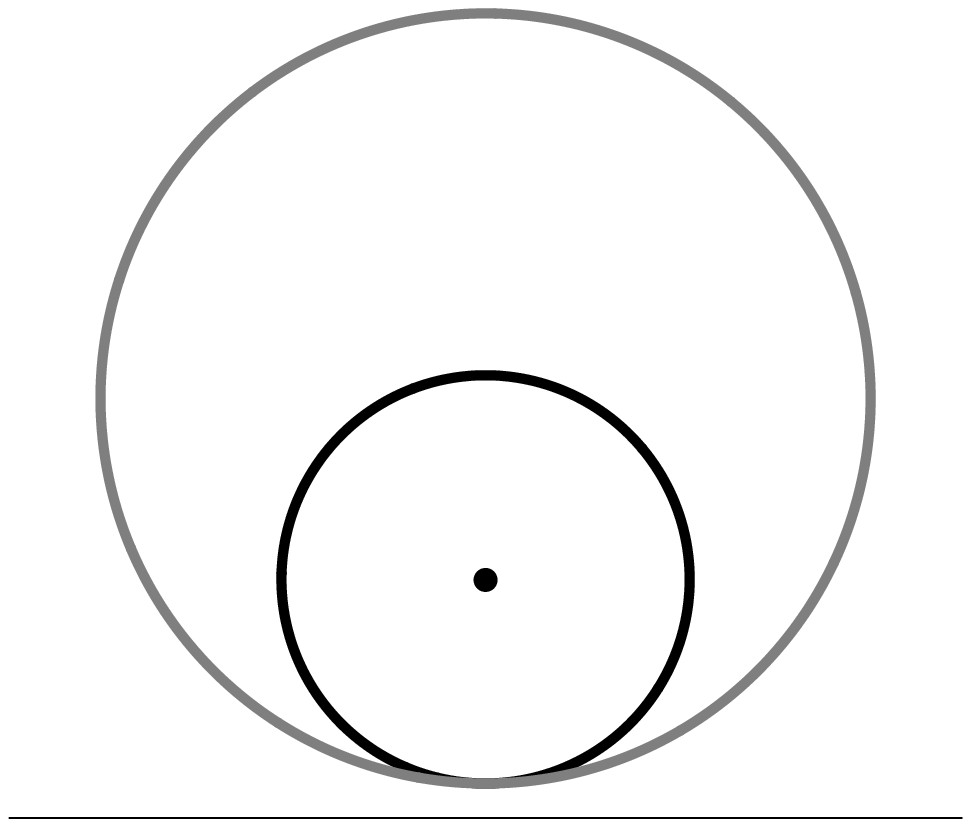}
    \caption{An example of inclusions of chordal disks (black) and quasihyperbolic disks (gray) in half-plane. The black line represents $\partial G$. The radius of the chordal disks is 0.45. \label{fig5}}
  \end{center}
\end{figure}

Results of Theorem \ref{kqk2sharp} and Corollary \ref{qkq2sharp} are illustrated in the planar case in Figure \ref{fig5}.

Radii $m$ and $M$ of Theorem \ref{jqj2}, Corollary \ref{qjq2}, Theorem \ref{kqk2sharp} and Corollary \ref{qkq2sharp} depend on $r$ and $|x|$. We can use \eqref{boundsforx} to make some of the radii independent of $|x|$.

\begin{theorem}\label{main2}
  Let $G = \Hn$, $x \in G$ with $x_1 = x_2 = \cdots = x_{n-1} = 0$ and $r \in (0,r_q(x))$. Then
  \[
    B_j(x,m_1) \subset B_q(x,r) \quad \textrm{and} \quad B_k(x,m_2) \subset B_q(x,r),
  \]
  where
  \[
    m_1 = \log \left( 1+ \frac{2r}{1-r^2} \right) \quad \textrm{and} \quad m_2 = \log \left( 1+\frac{2r^2}{\sqrt{1-r^2}} \right).
  \]
\end{theorem}
\begin{proof}
  Radius $m_2$ follows from Theorem \ref{main1}.

  By Theorem \ref{jqj2} we need to show that $m_1 \le \min \{ f(t),g(t) \}$ for $t \in (r/\sqrt{1-r^2},\sqrt{1-r^2}/r)$ and functions
  \[
    f(t) = \log \left( 1+\frac{r (1+ t^2)}{t\sqrt{1-r^2-r^2 t^2}} \right),\quad g(t) = \log \left( 1+\frac{r(1+t^2)}{\sqrt{1-r^2}t-r} \right).
  \]
  By elementary computation functions $f'(t)$ and $g'(t)$ are monotone and $f'(t_f)=0$ is equivalent to $t_f=\sqrt{1-r^2}/\sqrt{1+r^2}$ and $g'(t_g)=0$ is equivalent to $t_g=(1+r)/\sqrt{1-r^2}$. Now
  \[
    g(t_g) = \log \left( 1+\frac{2r}{1-r} \right) \ge \log \left( 1+\frac{2r}{1-r^2} \right) = f(t_f)
  \]
  and the assertion follows.
\end{proof}

\begin{corollary}
  Let $G = \Hn$, $x \in G$ with $x_1 = x_2 = \cdots = x_{n-1} = 0$ and $r > 0$. Then
  \[
    B_j(x,r) \subset B_q(x,M_1)
  \]
  and for $r \in (0,\log(1+\sqrt{2}))$
  \[
    \quad B_k(x,r) \subset B_q(x,M_2),
  \]
  where
  \[
    M_1 = \frac{\sqrt{2+e^r(e^r-2)}-1}{e^r-1}
  \]
  and
  \[
    M_2 = \frac{\sqrt{(e^r-1)(e^r-1+\sqrt{17+e^r(e^r-2)})}}{2\sqrt{2}}.
  \]
\end{corollary}
Note that for $M$ such that $B_q(x,r) \subset B_j(x,M)$ we have by Theorem \ref{jqj2} that $M \to \infty$ whenever $|x| \to r/\sqrt{1-r^2}$ or $|x| \to \sqrt{1-r^2}/r$. Similarly by Theorem \ref{kqk2sharp} we have that $M \to \infty$ whenever $|x| \to r/\sqrt{1-r^2}$ or $|x| \to \sqrt{1-r^2}/r$ for $M$ such that $B_q(x,r) \subset B_j(x,M)$.

\begin{proof}[Proof of Theorem \ref{mainthm}]
  The assertion follows from Theorems \ref{jkj}, \ref{main1}, \ref{jkj2} and \ref{main2}.
\end{proof}

\textbf{Acknowledgements.} The first author gratefully acknowledges the financial support of the Academy of Finland project 120972 of professor Pekka Koskela.



\end{document}